\documentclass[a4paper,11pt]{article}

\usepackage[utf8]{inputenc}
\usepackage[english,activeacute]{babel}
\usepackage{amsmath,amsthm,amsfonts,amssymb}
\usepackage{mathpartir}
\usepackage{ stmaryrd }
\usepackage{graphics}
\usepackage{enumerate}
\usepackage{mathrsfs}
\usepackage{rotating}

\input{xy}
\usepackage[all]{xy}
\usepackage{tikz}
\usepackage{tikz-cd}

\newtheorem{thm}{Theorem}[section]
\newtheorem{cor}[thm]{Corollary}

\theoremstyle{plain}

\theoremstyle{remark}

\input{diagxy}
\xyoption{curve}

\newbox\anglebox 
\setbox\anglebox=\hbox{\xy \POS(75,0)\ar@{-} (0,0) \ar@{-} (75,75)\endxy}
 
\newbox\angleboxr 
\setbox\angleboxr=\hbox{\xy \POS(0,0)\ar@{-} (0,75) \ar@{-} (75,0)\endxy}
 
\newbox\sanglebox 
\setbox\sanglebox=\hbox{\xy \POS(50,0)\ar@{-} (0,0) \ar@{-} (50,50)\endxy}
 
\newbox\sangleboxr 
\setbox\sangleboxr=\hbox{\xy \POS(0,0)\ar@{-} (0,50) \ar@{-} (50,0)\endxy}
 
\newbox\sangleboxf 
\setbox\sangleboxf=\hbox{\xy \POS(50,50)\ar@{-} (50,0) \ar@{-} (0,50)\endxy}
 
\newbox\angleboxf 
\setbox\angleboxf=\hbox{\xy \POS(75,75)\ar@{-} (75,0) \ar@{-} (0,75)\endxy}
 
\newbox\sangleboxfr 
\setbox\sangleboxfr=\hbox{\xy \POS(0,50)\ar@{-} (50,50) \ar@{-} (0,0)\endxy}
 
\newbox\angleboxfr 
\setbox\angleboxfr=\hbox{\xy \POS(0,75)\ar@{-} (75,75) \ar@{-} (0,0)\endxy}

\newcommand{\Sets}{\ensuremath{\mathbf{Set}}}

\title{A short proof of Shelah's eventual categoricity conjecture for AEC's with interpolation, under $GCH$}
\author{Christian Esp\'indola}

\begin{document}
\date{}
\maketitle

\begin{abstract}
We provide a short proof of Shelah's eventual categoricity conjecture, assuming the Generalized Continuum Hypothesis ($GCH$), for abstract elementary classes (AEC's) with interpolation, a strengthening of amalgamation which is a necessary and sufficient condition for an AEC categorical in a high enough cardinal to satisfy eventual categoricity. The proof builds on an earlier topos-theoretic argument which was syntactic in nature and recurred to $\kappa$-classifying toposes. We carry out here the same proof idea but from the semantic perspective, making use of a connection between $\kappa$-classifying toposes on one hand and the Scott adjunction on the other hand, this latter developed independently.\footnote{This research has been supported through the grant 19-00902S from the Grant Agency of the Czech Republic. I would like to thank Ji\v{r}\'i Rosick\'y and Ivan Di Liberti for reading and commenting a first version.}
\end{abstract}

\section{Introduction}
  Shelah's eventual categoricity conjecture asserts that for any AEC there is a cardinal $\kappa$ such that if the AEC is categorical in some $\lambda>\kappa$, it is categorical in \emph{all} $\lambda>\kappa$. This general conjecture was stated in \cite{shelah}, while the version for the particular case of sentences in $\mathcal{L}_{\omega_1, \omega}$ was conjectured circa 1977. Both conjectures are still open, though several approximations are known (see, e.g., \cite{sv3}). When the AEC has interpolation, the conjecture was proven to be true in \cite{espindolas} through a topos-theoretic argument, and in fact interpolation was shown to be a necessary and sufficient condition for eventual categoricity. We will here run the same proof, but taking a look at the semantic content of the $\kappa$-classifying toposes there employed. This is possible through an unexpected connection between these toposes and the ones arising from the Scott adjunction, this latter discovered by Simon Henry and Ivan Di Liberti (see \cite{sh} and Di Liberti's PhD thesis \cite{idl}). The use of the semantic counterpart of the $\kappa$-classifying toposes allows us to shed more light on the proof.

  The interpolation condition states that for an arbitrary set $\mathbb{S}$ of $\kappa$-coherent sentences over a language $\mathcal{L} \cup \mathbf{c}$ and any $\kappa$-coherent sentence $\phi$ over the language $\mathcal{L} \cup \mathbf{d}$ (where $\mathbf{c}, \mathbf{d}$ are disjoint sets of fresh constants), if $\mathbb{S} \vDash_{\geq \kappa} \phi$ (i.e., every model of size at least $\kappa$ of the accessible category with an expansion to a model of $\mathbb{S}$ is a model of $\phi$ in every expansion), then there is a $\kappa$-coherent interpolant $\theta$ over $\mathcal{L}$, i.e., such that $\mathbb{S} \vDash_{\geq \kappa} \theta \vDash_{\geq \kappa} \phi$.

  We remark that the use of $GCH$ made in this note is only meant for expository purposes as a simplification of the proofs,  but can be avoided throughout through the same forcing arguments of \cite{espindolas}.

\section{The $\lambda$-Scott topos}

The theorem that allows the topos-theoretic machinery to work is the following:

\begin{thm}\label{main}
Assume $\kappa^{<\kappa}=\kappa$. Then any $\kappa$-separable $\kappa$-topos has enough $\kappa$-points.
\end{thm}

For a proof, see \cite{espindolad}. Using this theorem, we can now find an important link between categoricity and atomic toposes:

\begin{thm}\label{cat}
Assume $\kappa^{<\kappa}=\kappa$. Any two $\kappa$-points of size at most $\kappa$ of a $\kappa$-separable $\kappa$-topos are $\mathcal{L}_{\infty, \kappa}$-elementarily equivalent if and only if the topos is two-valued and Boolean. 
\end{thm}

\begin{proof}
The proof is a direct generalization of the case $\kappa=\omega$, proven by Barr and Makkai in \cite{bm}, using Theorem \ref{main} to generalize the analogous statement for separable toposes (which is the particular case $\kappa=\omega$). We refer to \cite{espindolas} for details.
\end{proof}

Recall from \cite{sh} and \cite{idl} that there is an adjunction between $\kappa$-exact localizations of presheaf toposes and accessible categories with $\kappa$-directed colimits. Given such an accessible category $A$, its $\kappa$-Scott topos $S_{\kappa}(A)$ is the topos of functors to $\mathcal{S}et$ preserving $\kappa$-directed colimits. This assignment is functorial and gives one part of the adjunction, the other being the functor $pt_{\kappa}$ sending a topos to its category of $\kappa$-points. 

  Any AEC $\mathcal{K}$ is a $\lambda$-accessible category with directed colimits, and thus we can consider the corresponding Scott topos. Every $AEC$ is isomorphic to the category of models of the $\mathcal{L}_{\kappa, \kappa}$-theory axiomatizing its substructure functorial expansion (plus axioms defining an inequality predicate and complements of atomic formulas), and moreover, internal size coincides with cardinality at any $\lambda$.

\section{Classifying toposes for saturated models}

  The following theorem and proof appear in \cite{espindolas}; we reproduce it here for the sake of completeness:

\begin{thm}\label{sat}
Assume $GCH$ and interpolation. Then the $\kappa^+$-classifying topos $\Sets[\mathbb{T}^{sat}_{\kappa^+}]_{\kappa^+}$ of the theory of $\kappa^+$-saturated models is precisely $\mathcal{S}h(\mathcal{K}_{\kappa}^{op}, \tau_D)$, where $\tau_D$ is the dense topology (i.e., it is the double negation subtopos of $\Sets^{\mathcal{K}_{\kappa}}$).
\end{thm}

\begin{proof}
Consider the following diagram:

\begin{center}
\begin{tikzcd}
(\mathcal{C}_{\mathbb{T}_{\kappa}})_{\kappa^+} \arrow[rr, "ev"]  &  & \Sets^{\mathcal{K}_{\kappa}} \arrow[dd, "f^*"', bend right] \arrow[rrd, "M \cong \varinjlim_iev_{N_i}"] &  &               \\
  &  &   &  & \Sets \\
  &  & {\mathcal{S}h(K_{\kappa}^{op}, \tau_D)} \arrow[rru, dashed] \arrow[uu, "f_*"', bend right]                        &  &              
\end{tikzcd}
\end{center}
\noindent
Since we can assume that we have amalgamation (which is an easy consequence of interpolation), we have that a model $M$ is $\kappa^+$-saturated if and only if for every $p: N \to N'$ in $\mathcal{K}_{\kappa}$, each $q: N \to M$ extends to some $q': N' \to M$:
\begin{center}
\begin{tikzcd}
N' \arrow[rrd, "q'", dashed] &  &   \\
N \arrow[rr, "q"',] \arrow[u, "p",] &  & M
\end{tikzcd}
\end{center}
\noindent
This is the same as saying that $M: \Sets^{\mathcal{K}_{\kappa}} \to \Sets$ maps $p^*: [N', -] \to [N, -]$ to an epimorphism, since, writing $M \cong \varinjlim_iev_{N_i}$ as a $\kappa^+$-filtered colimit of evaluations, and given that any $N$ in $\mathcal{K}_{\kappa}$ is $\kappa^+$-presentable, we have:

$$\varinjlim_iev_{N_i}([N, -]) = \varinjlim_i [N, N_i] \cong [N, \varinjlim_i N_i] \cong [N, M]$$
\noindent
It follows that $M$ is $\kappa^+$-saturated if and only if it factors through $f^*: \Sets^{\mathcal{K}_{\kappa}} \to \mathcal{S}h(\mathcal{K}_{\kappa}^{op}, \tau_D)$, since the amalgamation property in $\mathcal{K}_{\kappa}$ implies that the dense topology coincides with the atomic topology (every nonempty sieve covers). This finishes the proof.
\end{proof}

\section{Eventual categoricity}

We are now going to prove the following:
  
\begin{thm}\label{btwn}
Assume $GCH$. Suppose the AEC $(\mathcal{K}, \prec)$ has interpolation and is categorical in some pair of cardinals $\lambda > \kappa > \theta$. Then $\mathcal{K}$ is also categorical at any $\delta$ with $\lambda \geq \delta \geq \kappa$.
\end{thm}

\begin{proof}
Take first $\delta=\kappa^+$, and let us see that $\Sets[\mathbb{T}_{\kappa^+}]_{\kappa^+}$ has a unique $\kappa^+$-point of size $\kappa^+$, up to isomorphism. Assume first that $\lambda$ is regular and consider the following diagram of toposes and inverse images of geometric morphisms given by restriction:

\begin{center}
\begin{tikzcd}
S_{\kappa^+}(\mathcal{K}_{\geq \kappa^+}) \arrow[rrr] \arrow[dd] &  &  & S_{\lambda}(\mathcal{K}_{\geq \lambda}) \arrow[dd, "\cong"] \\
                                                                 &  &  &                                                             \\
S_{\kappa^+}(Sat_{\kappa^+}(\mathcal{K})) \arrow[rrr]            &  &  & S_{\lambda}(Sat_{\lambda}(\mathcal{K}))                    
\end{tikzcd}
\end{center}

Let $Sat_{\kappa^+}(\mathcal{K})$ be axiomatized by the theory $\mathbb{T}^{sat}_{\kappa^+}$. Then the morphism $\mathcal{C}_{\mathbb{T}_{\kappa^+}} \to \mathcal{C}_{\mathbb{T}^{sat}_{\kappa^+}}$ induces a morphism between the corresponding $\kappa^+$-classifying toposes $f^*: \mathcal{S}et[\mathbb{T}_{\kappa^+}]_{\kappa^+} \to \mathcal{S}et[\mathbb{T}^{sat}_{\kappa^+}]_{\kappa^+}$. Then:

\begin{center}
\begin{tikzcd}
{ \mathcal{S}et[\mathbb{T}_{\kappa^+}]_{\kappa^+}} \arrow[dd, "f^*"'] \arrow[rr, "{\eta^*_{\mathcal{S}et[\mathbb{T}_{\kappa^+}]_{\kappa^+}}}"] &  & S_{\kappa^+}(\mathcal{K}_{\geq \kappa^+}) \arrow[rrr] \arrow[dd] &  &  & S_{\lambda}(\mathcal{K}_{\geq \lambda}) \arrow[dd, "\cong"] \\
                                                                                                                                               &  &                                                                  &  &  &                                                             \\
{\mathcal{S}et[\mathbb{T}^{sat}_{\kappa^+}]_{\kappa^+}} \arrow[rr, "{\eta^*_{\mathcal{S}et[\mathbb{T}^{sat}_{\kappa^+}]_{\kappa^+}}}"]         &  & S_{\kappa^+}(Sat_{\kappa^+}(\mathcal{K})) \arrow[rrr]            &  &  & S_{\lambda}(Sat_{\lambda}(\mathcal{K}))                    
\end{tikzcd}
\end{center}
By the results of \cite{rosicky}, $\mathcal{K}_{\geq \lambda}$ coincides with $Sat_{\lambda}(\mathcal{K})$. Therefore the right morphism is an isomorphism. We now want to deduce from this that $f^*$ is an equivalence, i.e., every model of size $\kappa^+$ is $\kappa^+$-saturated. Then $\mathcal{K}$ is $\kappa^+$-categorical since there is a unique such model (see \cite{rosicky}).

  We will prove that any sequent $\phi \vdash_{\mathbf{x}} \psi$ valid in $\mathcal{S}h(\mathcal{K}_{\kappa}^{op}, \tau_D)$ is already valid in $\Sets[\mathbb{T}_{\kappa^+}]_{\kappa^+}$. Note that such a sequent can be written, without loss of generality, as $\forall \mathbf{x} \phi(\mathbf{x})$ for $\kappa^+$-coherent $\phi$, and that it must be valid in the model of size $\lambda$. This is the case even if $\lambda$ is singular, as then categoricity implies $\kappa^+$-saturation by Lemma 4.3 in \cite{espindolas}.

  Now by the interpolation condition, whenever we have $\exists \mathbf{y} \bigwedge_{i \in S} \phi_i \vdash_{\mathbf{x}} \phi$, we can find some $\kappa^+$-coherent formula $\theta$ in $\mathcal{L}_{\kappa^{++}, \kappa^+}$ such that $\exists \mathbf{y} \bigwedge_{i \in S} \phi_i \vdash_{\mathbf{x}} \theta \vdash_{\mathbf{x}} \phi$ (indeed, it is equivalent to finding an interpolant of $\bigwedge_{i \in S} \phi_i(\mathbf{c}) \vdash_{\geq \kappa^+} \phi(\mathbf{d})$. Taking $\exists \mathbf{y} \bigwedge_{i \in S} \phi_i$ the formula obtained from the diagram of $N_{\lambda}$, we get the following strengthening of amalgamation: any morphism between models of size $\kappa$ $M \to M$ can be amalgamated with $M \to N_{\lambda}$ and additionally, if there is an embedding of $M$ into a model of bigger size satisfying $\neg \phi(\mathbf{e})$, we can take the amalgamation in such a way that the resulting model also satisfies $\neg \phi(\mathbf{e})$. It follows immediately from this that any model of size $\kappa^+$ must satisfy $\forall \mathbf{x} \phi$ and every model of size $\kappa^+$ is $\kappa^+$-saturated, as we wanted to prove. In the case that $\lambda$ is singular, we deduce from Lemma 4.3 in \cite{espindolas} that $N_{\lambda}$ is $\kappa^+$-closed, which is all we need to apply the above strengthening of the amalgamation property.

  Finally, note that the same proof above which allowed us to conclude categoricity in $\kappa^+$ from categoricity in $\kappa$, allows us now to conclude categoricity in $\kappa^{++}$, and so on. Categoricity in a $\delta$ which is a limit cardinal is easily handled knowing that $\mathcal{K}$ will be categorical at all $\gamma_i$ for a cofinal sequence of successors $\kappa<\gamma_i<\delta$. Indeed, since $\mathcal{K}$ is $\gamma_i$-categorical, the model of size $\gamma_i$ is $\gamma_i$-saturated, which allows us to successively find a set of compatible isomorphisms between submodels of any two models of size $\delta$ (using directed colimits at limit steps), proving that they are indeed isomorphic. 
\end{proof}

\begin{cor}\label{aec}
(Shelah's eventual categoricity conjecture for AEC's with interpolation). Assume $GCH$. Let $\mathcal{K}$ be an AEC with interpolation. Then there exists a cardinal $\mu_0$ such that if $\mathcal{K}$ is categorical in some $\lambda > \mu_0$, it is categorical in all $\lambda' > \mu_0$. Moreover, eventual categoricity for an AEC categorical in a high enough cardinal is equivalent to the AEC eventually satisfying interpolation.
\end{cor}

\begin{proof} Every $AEC$ is axiomatizable, and internal size coincides with cardinality at any $\lambda$. Our result follows taking $\mu_0$ to be the maximum of the Hanf numbers for categoricity and non-categoricity. When eventual categoricity holds, interpolation must hold eventually as seen in the proof of Theorem \ref{btwn}. 
\end{proof}

\bibliographystyle{amsalpha}

\renewcommand{\bibname}{References} 

\bibliography{references}



\end{document}